\documentclass{article}

\usepackage[utf8]{inputenc}
\usepackage[T1]{fontenc}
\usepackage{graphicx} 
\usepackage{amsmath,mathtools,amsfonts,amsthm}
\usepackage{siunitx}
\usepackage{tikz}
\usepackage[pdftex,colorlinks=true,linkcolor=blue,citecolor=blue,urlcolor=blue]{hyperref}

\newcommand\R{\mathbb{R}}
\newcommand\D{\mathrm{d}}

\newcommand{\ComR}[1]{{\color{black} #1}}

\newcommand{\Comus}[1]{{\color{black} #1}}

\newtheorem{theorem}{Theorem}[section]
\newtheorem{lemma}{Lemma}[section]
\newtheorem{remark}{Remark}[section]

\title{Topography optimization for enhancing microalgal growth in raceway ponds}
                                                
\author{Olivier Bernard$^{1}$, Liu-Di Lu$^{2}$, Jacques Sainte-Marie$^{3}$, \\
Julien Salomon$^{3}$}

\date{%
\noindent{\small\textit{$^1$NRIA Sophia Antipolis M\'editerran\'ee, BIOCORE Project-Team, Universit\'e Nice C\^ote d'Azur, 2004, Route des Lucioles - BP93, 06902 Sophia-Antipolis Cedex, France and Sorbonne Universit\'e, INSU-CNRS, Laboratoire d’Oc\'eanographie de Villefranche, 181 Chemin du Lazaret, 06230 Villefranche-sur-mer, France}}\\
{\small\textit{$^2$Section de Mathématiques, Université de Genève, rue du Conseil-Général 5-7, CP 64, 1205, Geneva, Switzerland}}\\
{\small\textit{$^3$INRIA Paris, ANGE Project-Team, 75589 Paris Cedex 12, France and Sorbonne Universit\'e, CNRS, Laboratoire Jacques-Louis Lions, 75005 Paris, France}}\\%
}

\begin{document}

\maketitle

\begin{abstract}
Modelling the evolution process for the growth of microalgae in an artificial pond is a huge challenge, given the complex interaction between hydrodynamics and biological processes occurring across various timescales. In this paper, we consider a raceway, i.e., an oval pond where the water is set in motion by a paddle wheel. Our aim is to investigate theoretically and numerically the impact of bottom topography in such raceway ponds on microalgae growth. To achieve this goal, we consider a biological model based on the Han model, coupled with the Saint--Venant systems that model the fluid. We then formulate an optimization problem, for which we apply the weak maximum principle to characterize optimal topographies that maximize biomass production over one lap of the raceway pond or multiple laps with a paddle wheel. In contrast to a widespread belief in the field of microalgae, we show that a flat topography in a periodic regime satisfies the necessary optimality condition, and observe in the numerical experiments that the flat topography is actually optimal in this case. However, non-trivial topographies may be more advantageous in alternative scenarios, such as when considering the effects of mixing devices within the model. This study sheds light on the intricate relationship between bottom topography, fluid dynamics, and microalgae growth in raceway ponds, offering valuable insights into optimizing biomass production.
\end{abstract}

\begin{paragraph}{Keyword: }
optimal control, weak maximum principle, microalgae, Han model, Saint--Venant system, raceway pond, shape optimization  
\end{paragraph}

\section{Introduction}

The numerical design of microalgae production technologies has been for decades a source of many interesting challenges not only in engineering but also in the area of scientific computing~\cite{Eilers1993, Masojidek2003, Yoo2012, Lamare2019}. The potential of these emerging photosynthetic organisms is found in cosmetics,  pharmaceutical fields, food, and - in the long term - in green chemistry and energy applications~\cite{Wijffels2010}. Outdoor production is mainly carried out in open bioreactors with a raceway shape. Algae grow while exposed to solar radiation in these circular basins, where the water is set in motion by a paddle wheel. This mixing device homogenizes the medium,  ensures equidistribution of nutrients, and guarantees that each cell will have regular access to light~\cite{Chiaramonti2013, Demory2018}. The algae are harvested periodically, and their concentration is maintained around an optimal value~\cite{Munoz-Tamayo2013, Posten2016}. The penetration of light is strongly reduced by the algal biomass, and less than 1\% of the incident light reaches the reactor bottom~\cite{Bernard2015}. In the case of larger biomass,  the light extinction is so high that a large fraction of the population evolves in the dark and does not grow anymore. At low biomass density,  a fraction of the solar light is not used by the algae and the productivity is suboptimal. Theoretical work has determined the optimal biomass for maximizing productivity~\cite{Masci2010, Grognard2014, Bernard2021optimal}.

Here, we consider another approach which consists in improving the photoproduction process by controlling the cell trajectories in the light field. We start from the observation that algal raceway ponds are dynamical systems combining a physical aspect - the hydrodynamical behavior of the fluid transporting the algae culture, and a biological aspect - the light harvesting by chlorophyll complexes in the cells~\cite{Bernard2013, Olivieri2015, Papacek2018}. We then study the effect of topography (or bathymetry) on growth to optimize the light received by the microalgae. Modelling this system is challenging, since it also involves the free-surface incompressible Navier--Stokes system~\cite{Bouharguane2012, Cocquet2021, Dongeren2008, Mohammadi2011}. The complexity of this model generally prevents \Comus{from} obtaining explicit formulas, and large computational resources are required to perform simulations.

Several experimental campaigns~\cite{Mendoza2013, Prussi2014} have shown that in straight sections \Comus{of} the raceway, the flow is not disturbed (which was further confirmed by CFD modelling~\cite{Hreiz2014, Inostroza2021}). Therefore, in these regions, despite turbulent dispersion, mixing is relatively poor. This mixing is mainly induced locally by the paddle wheel and, to a lesser extent, by the bends. The recent study of~\cite{Inostroza2021} confirms this finding, i.e.,  the turbulence is mainly generated near the paddle wheel and close to the surface.

We therefore focus on the main part of the raceway, outside the paddle-wheel area, and assume laminar flux. We study how to improve productivity in this part by modifying the bottom topography. This enables us to discuss the common belief that some specific topographies can bring more light to the algae in lower parts of the raceway, since cells get closer to the surface when reaching peaks in these topographies.

Let us detail our approach. We first introduce a coupled model to represent the growth of algae in a one-dimensional (1D) raceway pond, accounting for the light that they receive. This model is obtained by combining the Han photosynthesis equations with a hydrodynamic law based on the Saint--Venant system. This first step enables us to formulate an optimization problem in which the topography of the raceway is designed to maximize productivity. We then use an adjoint-based optimization scheme to include the constraints associated with the Saint--Venant regime. We prove that the flat topography satisfies the first-order optimality systems in a periodic case, focusing on the fraction of the raceway in laminar regime. However, non-trivial topographies can be obtained in other contexts, e.g., when the periodic assumption is removed or when the mixing device is accounted for in conjunction with the bottom topography. Numerical simulations show that a combination of turbulence-induced mixing and non-flat topographies can slightly increase biomass production. However, enhancing the turbulence by mixing significantly increases productivity and is definitely the most efficient approach~\cite{Bernard03170481, Bernard02970756}, even if more energy is dissipated in this process.

The outline of the paper is as follows. In Section~\ref{sec:model},  we present the biological and hydrodynamical models underlying our coupled system. In Section~\ref{sec:opti}, we describe the optimization problem and a corresponding numerical optimization procedure. Section~\ref{sec:num} is devoted to the numerical results obtained with our approach. We then conclude with some perspectives opened up by this work.

\section{Hydrodynamic and biological models}\label{sec:model}

Our approach is based on a coupling of the hydrodynamic transport of the particles with the photosystem evolution driven by the light intensity they receive when traveling in the raceway pond.

\subsection{Hydrodynamical model and Lagrangian trajectories}\label{sec:model_hydro}

Saint--Venant equations are a popular model of geophysical flows. This system is derived from the free surface incompressible Navier--Stokes equations (see, for instance,~\cite{Gerbeau2001}). Here, we focus on its 1D smooth steady state solutions in a laminar regime, which satisfy
\begin{equation}\label{eq:sv}
\partial_x(h u) = 0,  \quad \partial_x(h u^2+g\frac{h^2}2) = -g h \partial_x z_b,
\end{equation}
where $h$ is the water depth, $u$ is the horizontal averaged velocity of the fluid, the $g$ is the gravitational constant, and $z_b$ is the topography. The free surface $\eta$ and the average discharge are given by $\eta := h + z_b$ and $Q = h u$ respectively. This system is presented in Figure~\ref{fig:shallow}.
\begin{figure}
\begin{center}
\begin{tikzpicture}[scale=.8]
\draw [thick](0, -4) -- (6, -4);
\draw (0, -4) node[anchor=north] {0};
\draw (6, -4) node[anchor=north] {$L$};
\draw[thick,  ->] (0, -4) -- (0, 0.5) node[anchor=east] {$z$};
\draw (0, -0.2) node[anchor=east] {$0$};
\draw[dashed, ->] (0, -0.2) -- (6.2, -0.2) node[anchor=west]  {$x$};
\draw [dashed,  <-](1.5, -0.88) -- (1.5, -0.2);
\draw (1.5, 0.2) node{$\eta(x)$};
\draw[thick, <-] (2.5, -0.65) -- (2.5, 0.5);
\draw[thick, <-] (3, -0.4) -- (3, 0.5);
\draw[thick, <-] (3.5, -0.2) -- (3.5, 0.5);
\draw[thick, <-] (4.0, -0.08) -- (4, 0.5);
\draw[thick, <-] (4.5, -0.03) -- (4.5, 0.5);
\draw[thick, <-] (5, -0.1) -- (5, 0.5);
\draw[thick, <-] (5.5, -0.28) -- (5.5, 0.5);
\draw (4, 0.7) node {$I_s$};
\draw (0, -0.4) sin (1.5, -0.9) cos (3, -0.4) sin (4.5, 0.0) cos (6, -0.6);
\draw [dashed] (0, -1.4) sin (1.5, -1.9) cos (3, -1.4) sin (4.5, -1) cos (6, -1.6);
\draw [dashed] (0, -2.4) sin (1.5, -2.9) cos (3, -2.4) sin (4.5, -2) cos (6, -2.6);
\draw[fill,  gray] (0, -3.1) sin (1.5, -3.8) cos (3, -3.3) sin (4.5, -2.7) cos (6, -3.3) -- (6, -4) -- (0, -4) ;
\draw [dashed, ->](5, -0.2) -- (5, -2.75);
\draw (5, -1.8) node[anchor=west] {$z_b(x)$};
\draw[dashed, <->] (3, -0.4) -- (3, -3.3);
\draw (3, -1.85) node[anchor=west]{$h(x)$};
\node at (2, -2.3){$u(x)$};
\draw [->](1.3, -1.4) -- (1.6, -1.4);
\draw [->](1.3, -2.4) -- (1.6, -2.4);
\draw [->](1.3, -3.4) -- (1.6, -3.4);
\end{tikzpicture}
\end{center}
\caption{Representation of the one dimensional hydrodynamic model.}
\label{fig:shallow}
\end{figure}
The $z$ (resp. $x$) axis represents the vertical (resp. horizontal) direction and $I_s$ is the light intensity on the free surface (assumed to be constant).

Integrating the equation on the left of~\eqref{eq:sv}, we get
\begin{equation}\label{eq:Q0}
h u = Q_0,  
\end{equation}
for a fixed positive constant $Q_0$. This implies a constant discharge in space. Then the equation on the right-hand side of~\eqref{eq:sv} can be rewritten by
\begin{equation}\label{eq:sv2mod}
h u\partial_x u + h\partial_x g h + h\partial_x g z_b = 0.
\end{equation}
Assume that $h$ is non-zero, dividing then the equality~\eqref{eq:sv2mod} by $h$ and using~\eqref{eq:Q0} to eliminate $u$, we get 
$\partial_x \Big( \frac{Q_0^2}{2h^2} + g (h + z_b)\Big) = 0.$ Given $h(0),  z_b(0)\in \R$, we obtain 
\[\frac{Q_0^2}{2h(x)^2} + g (h(x) + z_b(x)) = \frac{Q_0^2}{2h^2(0)} + g (h(0) + z_b(0)) =: M_0, \] 
which holds for all $x\in [0, L]$, meaning that the topography $z_b$ satisfies
\begin{equation}\label{eq:zb}
z_b = \frac{M_0}g - \frac{Q_0^2}{2g h^2} - h.
\end{equation}

\begin{remark}\label{rem:flu}
Let $Fr=\frac{u}{\sqrt{g h}}$ be the Froude number. The situation $Fr<1$ corresponds to the subcritical case (i.e., the flow regime is \emph{fluvial}), while $Fr>1$ corresponds to the supercritical case (i.e., the flow regime is \emph{torrential}). In the steady case, the threshold value $h=h_c$ is obtained for $Fr=1$; using~\eqref{eq:Q0}, we find $h_c := (\frac{Q^2_0}{g})^{\frac 13}$.
\end{remark}

Because of~\eqref{eq:zb}, $h$ solves a third-order polynomial equation. Given a smooth topography $z_b$, if $h_c + z_b + \frac{Q_0^2}{2g h^2_c} - \frac{M_0}g < 0,$ there exists a unique positive smooth solution of~\eqref{eq:zb} that satisfies the subcritical flow condition (see~\cite[Lemma 1]{Michel-Dansac2016}).

From the incompressibility of the flow, we have $\nabla \cdot \underline{\textbf{u}} = 0$ with $\underline{\textbf{u}} = (u(x), w(x, z))$. Here, $w(x, z)$ is the vertical velocity. Incompressibility implies $\partial_x u + \partial_z w = 0$. Integrating the latter from the topography $z_b$ to an arbitrary vertical position $z$ gives:
\[\begin{aligned}
0 & = \int_{z_b}^z \big(\partial_x u(x) + \partial_{\xi} w(x, \xi) \big)\D \xi \\
& =  (z - z_b)\partial_x u(x) + w(x, z) - w(x, z_b) \\
& = (z - z_b) \partial_x u(x) - u(x) \partial_x z_b + w(x, z)\\
& = \partial_x \big((z-z_b)u(x)\big) + w(x, z), 
\end{aligned}\]
where we have used the kinematic condition at the bottom, i.e., $w(x, z_b) = u(x)\partial_x z_b$. It follows from~\eqref{eq:zb} that
\begin{equation}\label{eq:w}
w(x , z) = \big(\frac{M_0}g - \frac{3u^2(x)}{2g} - z\big)u'(x), 
\end{equation}
with $u'(x)$ the derivative of $u$ with respect to $x$.

Let the pair $(x(t), z(t))$ be the position of a particle (or an algal cell) at time $t$ in the raceway pond. The Lagrangian trajectory is characterized by 
\begin{equation}\label{eq:xtzt}
\begin{pmatrix}
\dot{x}(t)\\
\dot{z}(t)   
\end{pmatrix}
= 
\begin{pmatrix}
u(x(t))\\
w\big(x(t), z(t)\big)  
\end{pmatrix}, 
\end{equation}
with the initial position at time $t = 0$, $(x(0), z(0)) = (x_0, z_0)$.

\begin{remark}
The geometry of the raceway pond with small dissipation and shear effects (reduced wall friction and viscosity) justifies a laminar flow modeled by a shallow-water model, such as the Saint--Venant system. This regime also minimizes the mixing energy and hence is favored at the industrial scale.

A higher mixing energy would lead to a turbulent regime. A possible way to enrich the representation of Lagrangian trajectories in this case would consist \Comus{of} including a Brownian in~\eqref{eq:xtzt}. However, getting time-free expressions of the trajectories (as in~\eqref{eq:z'x} and~\eqref{eq:c'x}) in this case is much more challenging, so that such a strategy would require a large set of simulations together with an averaging strategy.
\end{remark}

The Lagrangian trajectory given by~\eqref{eq:xtzt} is a general formulation, which still holds when we change the hydrodynamical model. In our setting, we can find a time-free formulation of the Lagrangian trajectory. More precisely, we denote by $z(x)$ the depth of a particle at position $x$. From~\eqref{eq:w} and~\eqref{eq:xtzt}, we get
\begin{equation}\label{eq:z'x}
z' := \frac{\dot{z}}{\dot{x}} = \big(\frac {M_0}g -\frac{3u^2}{2g} - z\big)\frac{u'}{u}.
\end{equation}
From~\eqref{eq:Q0}, \eqref{eq:zb} and the definition of the free surface $\eta$, we have 
\[\eta=h+z_b = \frac{M_0}g - \frac{u^2}{2g},\] 
which implies that $\eta'=-u u'/g$. Multiplying then~\eqref{eq:z'x} on both sides by $u$, and using the formulation of $\eta$ and $\eta'$, one finds 
\[z' u + z u'= \big(\eta - \frac{u^2}g\big)u' = \eta u' + \eta' u,\] 
which implies that $(u(z-\eta))'=0$. Using again the identity~\eqref{eq:Q0}, one obtains $\eta(x) - z(x) = \frac{h(x)}{h(0)}\big(\eta(0) - z(0)\big)$. This equation shows that given the initial water depth $h(0)$ and the initial free surface position $\eta(0)$, the distance between a trajectory $z$ (starting from the position $z(0)$) and the free surface $\eta$ depends only on the water depth $h$. On the other hand, the time-free formulation of the trajectory reads 
\begin{equation}\label{eq:zx}
z(x) =\eta(x) - \frac{h(x)}{h(0)}\big(\eta(0) - z(0) \big).
\end{equation}
We will further exploit the property of this formulation in Section~\ref{sec:opti}.

\begin{remark}
Since $Q_0$ is chosen to be positive, $h$ is necessarily positive. Moreover, if $z(0)$ belongs to $[z_b(0), \eta(0)]$, then $z(x)$ belongs to $[z_b(x), \eta(x)]$. In particular, choosing $z(0)=z_b(0)$ in~\eqref{eq:zx} and using~\eqref{eq:Q0} give $z(x) = z_b(x)$. In the same way, we find that $z(x)=\eta(x)$ when $z(0) = \eta(0)$. 
\end{remark}

\subsection{Modeling the photosystems dynamics}

To describe the dynamics of photosystems, we use here the Han model~\cite{Han2002}. This model is generally considered to characterize the photosynthetic process of these subunits as they harvest photons and transfer their energy to the cell to fix CO\textsubscript{2}.

\subsubsection{The Han model}

The Han model is a compartmental model in which the photosystems are described by three different states: open and ready to harvest a photon ($A$),  closed while processing the absorbed photon energy ($B$), or inhibited if several photons have been absorbed simultaneously ($C$). The relation of these three states are schematically presented in Fig.~\ref{fig:han}.
\usetikzlibrary{shapes.geometric}
\tikzset{block/.style = {draw, rectangle, minimum height = 2.5em, minimum width = 2.5em}}
\begin{figure}
\begin{center}
\begin{tikzpicture}
\node at (0, 0)[thick, circle, draw, scale=1.5] (A){$A$};
\node at (3, 0)[thick, circle, draw, scale=1.5] (B){$B$};
\node at (6, 0)[thick, circle, draw, scale=1.5] (C){$C$};
\node at (1.5, 0) [thick, block, draw, scale=0.8](sigI){$\sigma I$};
\node at (4.5, 0) [thick, block, draw, scale=0.8](kdsigI){$k_d\sigma I$};
\node at (1.5, 1.5) [thick, block, draw, scale=0.8](tau1){$\tau^{-1}$};
\node at (4.5, 1.5) [thick, block, draw, scale=0.8](kr){$k_r$};
\node at (0, 2)(Pho1){Photon $I$};
\node at (3, 2)(Pho2){Photon $I$};
\draw (A)--(sigI)[thick, ->](sigI)--(B);
\draw (B)--(kdsigI)[thick, ->](kdsigI)--(C);
\draw (C) edge[thick, bend right=35] (kr)(B) edge[thick, <-, bend left=35] (kr);
\draw (B) edge[thick, bend right=35] (tau1)(A) edge[thick, <-, bend left=35] (tau1);
\draw [thick, ->](Pho1)--(A);
\draw [thick, ->](Pho2)--(B);
\end{tikzpicture}
\end{center}
\caption{Han's model, describing the state transition probability, as a function of the photon flux.} 
\label{fig:han}
\end{figure}

The evolution satisfies the following ordinary differential equations (ODEs)
\begin{equation}\label{eq:Han}
\begin{aligned}
\dot A &= -\sigma I A + \frac B{\tau}, \\
\dot B &=  \sigma I A - \frac B{\tau} + k_rC - k_d\sigma I B, \\
\dot C &= -k_r C + k_d \sigma I B.
\end{aligned}
\end{equation}
Here, $I$ denotes the light density, a continuous time-varying signal. The states $A$, $B$, and $C$ are the relative frequencies of three possible states with $A+B+C=1$, so that~\eqref{eq:Han} can be reduced to a system in dimension two by eliminating the state $B$. Here, $\sigma$ stands for the specific photon absorption, $\tau$ is the turnover rate, $k_r$  and $k_d$ represent the photosystem repair and damage rates, which are all positive.

The dynamics of the open state $A$ can be shown to be much faster than the dynamics of the photoinhibition state $C$. A slow-fast approximation by using singular perturbation theory (as shown in details in~\cite{Lamare2019}) leads to the simplification of the dynamics driven by the slow dynamics of $C$:
\begin{equation}\label{eq:ct}
\dot C = -\alpha(I) C + \beta(I), 
\end{equation}
where 
\begin{equation}\label{eq:DefAlphaBeta}
\alpha(I) = k_d\tau \frac{(\sigma I)^2}{\tau \sigma I+1} + k_r, 
\quad 
\beta(I) =  k_d\tau \frac{(\sigma I)^2}{\tau \sigma I+1}.
\end{equation} 
Repeating the reasoning done to get~\eqref{eq:z'x} with~\eqref{eq:ct} and~\eqref{eq:Q0},  we can also find a time-free reformulation, namely
\begin{equation}\label{eq:c'x}
C':=\frac{\dot C}{\dot x} = \frac{-\alpha(I)C + \beta(I)}{Q_0}h, 
\end{equation}
where all the functions on the right-hand side only depend on the spatial variable $x$.

\subsubsection{Periodic setting} 
We consider the case where $C$ is periodic, with a period corresponding to one lap of the raceway pond. This situation occurs, e.g., when an appropriate harvest is performed after each lap. To describe the corresponding model, we first consider a variant of the usual Cauchy problem~\eqref{eq:c'x}:

{\it Given $I\in \mathcal C([0, L]; \;\R)$, $I\geq 0$,  find $(C_0, C)\in [0, 1]\times \mathcal C([0, L]; \;[0, 1])$ such that}
\begin{equation}\label{eq:Cperiodic}
\left\{
\begin{aligned}
C'(x) &= \frac{-\alpha\big(I(x)\big)C(x) + \beta\big(I(x)\big)}{Q_0}h(x),  \quad x\in[0, L], \\
C(L) &= C(0)= C_0.
\end{aligned}
\right.
\end{equation}
Let us show that the solution $C(x)$ of~\eqref{eq:Cperiodic} exists. Indeed, applying Duhamel's formula \Comus{to} the Cauchy problem associated with~\eqref{eq:c'x} and the initial condition $C(0) = C_0$, and using the inequality $\beta(I) \leq \alpha(I)$ gives
\[\begin{aligned}
C(L)-C_0 = &-\left(1-e^{-\int_0^{L} \frac{\alpha(I(s))h(s)}{Q_0}\, \D s}\right)C_0 + \int_0^{L} e^{-\int_s^{L} \frac{\alpha(I(y))h(y)}{Q_0}\, \D y} \frac{\beta\big(I(s)\big)h(s)}{Q_0}\, \D s\\
\leq & \left(1-e^{-\int_0^{L} \frac{\alpha(I(s))h(s)}{Q_0}\, \D s}\right)\left(1-C_0\right).
\end{aligned}\]
Hence, the affine mapping $\Phi: C_0\mapsto C(L)-C_0$ satisfies $\Phi(0) \geq 0$, and the inequality implies that $\Phi(1) \leq 0$. It follows that there exists a unique $C_0\in [0, 1]$ that satisfies $C(L) - C_0 = 0$. Using {\it Intermediate Value Theorem}, we get the next result.

\begin{theorem}
There exists a unique couple $(C_0, C) \in [0, 1]\times \mathcal C([0, L]; \; [0, 1])$ that satisfies~\eqref{eq:Cperiodic}.
\end{theorem}

\subsubsection{Growth rate}
Finally, the net growth rate of the photosystem is defined by balancing photosynthesis and respiration, which gives
\begin{equation}\label{eq:mu}
\mu(C, I) := \zeta(I)-\gamma(I)C, 
\end{equation}
where 
\begin{equation}\label{eq:DefZetaGamma}
\gamma(I) = \frac{k\sigma I}{\tau \sigma I+1}, 
\quad 
\zeta(I) = \frac{k\sigma I}{\tau \sigma I+1}-R.
\end{equation}
Here, $k$ is a factor that relates the received energy with the growth rate and $R$ represents the respiration rate.

\subsection{Coupling of two systems}

As shown in the previous section, the light intensity $I$ plays an important role in algal growth, since it triggers photosynthesis. On the other hand, the position of the algae influences the perceived light as well as the efficiency of the photosynthesis process. Therefore, light intensity is the main connection that couples the hydrodynamic model and the physiological evolution of algae. To evaluate the light intensity observed on the trajectory $z$, we assume that the growth process occurs at a much slower timescale than that of hydrodynamics and is, as such, negligible for one lap over the raceway. In the same way, uncertainties such as rainfall and evaporation, can also be neglected at this timescale. These factors can be taken into account for longer timescale using more detailed models, see for instance~\cite{DeLuca2017, Casagli2021}. In this framework, the Beer--Lambert law describes how light is attenuated with depth $\xi$ by $I(x, \xi) := I_s\exp\big(-\varepsilon (\eta(x) - \xi)\big)$, where $\varepsilon$ is the light extinction coefficient. Replacing $\xi$ in the previous formulation by the trajectory~\eqref{eq:zx}, we then get the following expression for the captured light intensity along the trajectory $z(x)$:
\begin{equation}\label{eq:Beer}
I\big(x, z(x)\big) = I_s \exp\Big(-\varepsilon  \frac{h(x)}{h(0)} \big( \eta(0) -z(0) \big) \Big).
\end{equation}
In particular, we observe that for given data $I_s$, $\varepsilon$, $h(0)$, and $\eta(0)$, the perceived light intensity along the trajectory $z(x)$ depends only on its initial position $z(0)$ and $h(x)$.

In order to evaluate the quality of this coupled system, we define the average net growth rate of the system by 
\begin{equation}\label{eq:mubar}
\bar \mu := \frac 1V\int_0^L\int_{z_b(x)}^{\eta(x)} \mu\big(C(x,  z),  I(x,  z)\big) \, \D z \D x, 
\end{equation}
where $\mu$ is defined by~\eqref{eq:mu} and $V := \int_0^L h(x)\D x$ is the volume of our 1D raceway.

\section{Optimization problem}\label{sec:opti}

In this section, we define the optimal control problems associated with our biological--hydrodynamic model. Depending on $V$, we divide our study into two cases.

\subsection{Objective function and vertical discretization}\label{sec:objfun}

Our goal is to find the optimal topography $z_b$ that maximizes the average net growth rate~\eqref{eq:mubar}. In order to tackle numerically this optimization problem, let us first consider a vertical discretization. Let $N_z$ denotes the number of trajectories, we consider a uniform vertical discretization of their initial position: 
\begin{equation}\label{eq:zi0}
z_i(0) := \eta(0) - \frac {i-\frac12}{N_z}h(0), 
\quad 
i=1, \ldots, N_z.
\end{equation}
Using the formulation~\eqref{eq:zx}, we find the trajectories $z_i(x) := \eta(x) - \frac {i-\frac12}{N_z}h(x)$, $i = 1, \ldots, N_z$. In particular, the distribution of \Comus{the} trajectories $z_i(x)$ remains uniform along the direction of $x$. \ComR{Using~\eqref{eq:zi0}}, we obtain the perceived light intensity on $z_i(x)$: 
\begin{equation}\label{eq:I_ih}
I\big(x, z_i(x)\big) = I_s \exp\Big(-\varepsilon \frac{h(x)}{h(0)} \big( \eta(0) -z_i(0) \big) \Big) 
= I_s \exp\Big(-\varepsilon \frac {i-\frac12}{N_z}h(x) \Big),
\end{equation}
where we use the closed form of the light intensity~\eqref{eq:Beer} and the definition of $z_i(0)$. To simplify \Comus{the} notation and emphasize the dependence on the water depth $h$, we write $I_i(h(x))$ instead of $I(x, z_i(x))$ hereafter. The photoinhibition state $C_i$ is then computed using the evolution~\eqref{eq:c'x} for $I = I_i(h)$. In this setting, the semi-discrete average net growth rate in the raceway pond can be derived from~\eqref{eq:mubar} as
\begin{equation}\label{eq:muNz}
\bar \mu_{N_z}\big(h\big) := \frac 1{V N_z}\sum_{i=1}^{N_z} \int_0^L\mu\Big( C_i(x), I_i\big(h(x)\big) \Big)h(x) \, \D x,
\end{equation}
where $h$ is the variable of the objective function, and $\mu$ is given by~\eqref{eq:mu}. From now on, we focus on the subcritical case, i.e., $Fr<1$, see Remark~\ref{rem:flu}. As mentioned in Section~\ref{sec:model_hydro}, in this regime, a given topography $z_b$ corresponds to a unique water depth $h$ which verifies this assumption.

\begin{remark}
Given a topography $z_b$, the usual shallow-water solvers typically consider equations of type~\eqref{eq:zb} to compute $h$ in the simulations. Here,  we use this equation in the opposite way, i.e., to recover $z_b$ from $h$. In this way, we directly optimize $h$ instead of $z_b$, since the expressions of the evolution of the state $C$~\eqref{eq:c'x}, the light intensity~\eqref{eq:I_ih} and the objective function~\eqref{eq:muNz} depend on $h$ and not on $z_b$.
\end{remark}

\subsection{Constant Volume}\label{sec:cst_vol}

For simplicity, we omit from now on the variable $x$ in the notation and consider $h$ as the variable of the light intensities $(I_i)_{i=1, \ldots, N_z}$ and $\bar \mu_{N_z}$. For a fixed volume $V>0$ and a discharge $Q_0>0$, we seek admissible controls $h \in L^\infty([0, L]; \, \R)$, $h>0$ over a fixed length $L>0$, which maximize the semi-discrete average net growth rate~\eqref{eq:muNz}. Thus, the optimal control problem (OCP) reads 
\begin{equation}\label{eq:ocp1}
\begin{aligned}
\max_{h\in L^\infty([0, L]; \; \R), \ h>0} \, \bar \mu_{N_z}(h) = \ &\sum_{i=1}^{N_z} \int_0^L \frac{\mu\Big(C_i(x), I_i\big(h(x)\big) \Big)}{V N_z}h(x) \, \D x,\\ 
\text{s.t.} \quad C'_i =\ & \frac{\beta\left(I_i(h)\right) - \alpha\left(I_i(h)\right)C_i}{Q_0}h, \\
C_i(0)=\ & C_i(L), 
\qquad 
\forall i=1, \cdots, N_z, \\
v'=\ & h, \\
v(0)=\ &0,  \ v(L)=V.
\end{aligned}
\tag{P1}
\end{equation}
Here, we use formula~\eqref{eq:mu} for $\mu$, $h$ is the control variable, and $(C_i, v)$ are the state variables, where $v$ has been introduced to take into account the constraint $V=hL$. The Hamiltonian associated with~\eqref{eq:ocp1} is given by
\[\begin{aligned}
H(C_i, v, p_{C_i}, p_v, p_0, h) = &\sum_{i=1}^{N_z} p_{C_i} \frac{\beta\left(I_i(h)\right) - \alpha\left(I_i(h)\right)C_i}{Q_0}h \\
&+ p_v h + p_0 \sum_{i=1}^{N_z} \frac{\ComR{\mu\Big(C_i, I_i\big(h\big) \Big)}}{V N_z}h, 
\end{aligned}\]
where $(p_{C_i}, p_v)$ are the co-states of $(C_i, v)$ respectively, and $p_0$ is a real number. Suppose that $h^{\star}\in L^\infty([0, L]; \, \R)$, \ComR{$h^{\star}>0$} is a maximizer, and $C_i^{\star}$, $v^{\star}$ are the corresponding solutions of the problem~\eqref{eq:ocp1}. Using the weak maximum principle~\cite[Pages 33--35]{Trelat2024}, there exist absolutely continuous functions $p_{C_i}^{\star}:[0, L]\to \R$, $p_v^{\star}:[0, L]\to \R$ and a real number $p_0^{\star}\leq 0$, such that for almost every $x\in[0, L]$, the extremals $(C_i^{\star}, v^{\star}, p_{C_i}^{\star}, p_v^{\star}, p_0^{\star}, h^{\star})$ satisfy the optimality system
\begin{equation}\label{eq:optsys}
\begin{aligned}
C_i' =\ & \frac{\partial H}{\partial p_{C_i}} 
= \frac{\beta\left(I_i(h)\right) - \alpha\left(I_i(h)\right)C_i}{Q_0}h, \quad 
v' = \frac{\partial H}{\partial p_v} = h,\\ 
{p_{C_i}}' =\ & -\frac{\partial H}{\partial C_i} 
= p_{C_i} \frac{\alpha\left(I_i(h)\right)}{Q_0}h + p_0\frac{\gamma\left(I_i(h)\right)}{VN_z}h,
\quad
{p_v}' = -\frac{\partial H}{\partial v} =0, \\ 
0 =\ & \frac{\partial H}{\partial h} = 
\sum_{i=1}^{N_z} p_{C_i}\frac{ \beta'\left(I_i(h)\right) - \alpha'\left(I_i(h)\right)C_i}{Q_0}I'_i(h)h 
+ \sum_{i=1}^{N_z} p_{C_i}\frac{\beta\left(I_i(h)\right) - \alpha\left(I_i(h)\right)C_i}{Q_0}\\
&+ p_0\sum_{i=1}^{N_z} \frac{\zeta'\left(I_i(h)\right) - \gamma'\left(I_i(h)\right)C_i}{VN_z}I'_i(h)h 
+ p_0\sum_{i=1}^{N_z} \frac{\zeta\left(I_i(h)\right) - \gamma\left(I_i(h)\right)C_i}{VN_z} + p_v. 
\end{aligned}
\end{equation}

\begin{lemma}\label{lem:normalVcst}
The extremal $(C_i^{\star}, v^{\star}, p_{C_i}^{\star}, p_v^{\star}, p_0^{\star}, h^{\star})$ that satisfies~\eqref{eq:optsys} is normal.
\end{lemma}

\begin{proof}
We use the equivalent dual form of the Mangasarian-Fromovitz constraint qualification~\cite[p. 255–269]{Solodov2011}, i.e., we prove that if $p_0^{\star} = 0$, then $p_{C_i}^{\star}$ and $p_v^{\star}$ are equal to zero on $[0, L]$.

Substituting $p_0^{\star} = 0$ into~\eqref{eq:optsys}, the ODE associated with $p_{C_i}^{\star}$ then reads
\begin{equation}\label{eq:adj}
({p_{C_i}^{\star}})' = p_{C_i}^{\star} \frac{\alpha\left(I_i(h^{\star})\right)}{Q_0}h^{\star}, 
\quad 
p_{C_i}^{\star}(0)= p_{C_i}^{\star}(L), 
\quad 
\forall i=1, \cdots, N_z,
\end{equation}
where we complete by the periodic condition determined using $C_i^{\star}(0)= C_i^{\star}(L)$, $\forall i=1, \cdots, N_z$. Note that $Q_0>0$ and $\alpha$ is a positive function from~\eqref{eq:DefAlphaBeta}, and $h^{\star}>0$. Hence, we have $\frac{\alpha\left(I_i(h^{\star})\right)}{Q_0}h^{\star} > 0$. Using then a similar reasoning to that for the system~\eqref{eq:Cperiodic}, we find that the only solution of~\eqref{eq:adj} is $p_{C_i}^{\star} = 0$. Substituting $p_{C_i}^{\star}=0$ and $p_0^{\star}=0$ into the last equation of~\eqref{eq:optsys}, we obtain $p_v^{\star}=0$, which contradicts the fact that $p_{C_i}^{\star}$ and $p_v^{\star}$ are not identically 0 on $[0, L]$. Therefore, $p_0^{\star}<0$.
\end{proof}

When the extremum is normal,  $p_{C_i}^{\star}$ and $p_v^{\star}$ are usually normalized so that $p_0^{\star} = -1$ \Comus{is} what we set hereafter. Let us show that the flat topography satisfies~\eqref{eq:optsys}.

\begin{theorem}\label{thm:flattopo}
There exists $p_v^f\in \R$ such that the constant water depth 
\[h^f:=\frac VL,\] 
and the corresponding solutions $(C_i^f)_{i=1, \cdots, N_z},  (p_{C_i}^f)_{i=1, \cdots, N_z}$, $v^f$ satisfy~\eqref{eq:optsys}.
\end{theorem}

\begin{proof}
From $v'=h^f$ with $v(0)=0$, $v(L)=V$, we find $v^f=\frac VL x$. Given $i\in\{1, \cdots, N_z\}$,  from~\eqref{eq:I_ih},  we deduce that 
\[I_i(h^f) = I_s \exp(-\varepsilon \frac{i-\frac 12}{N_z} h^f),  
\quad 
I_i'(h^f) = -\varepsilon \frac{i-\frac 12}{N_z} I_i(h^f),\] 
which are constant on $[0, L]$. Solving the equation of $C_i$ in~\eqref{eq:optsys} gives
\begin{equation}\label{eq:Cx}
\begin{aligned}
C_i(x) = e&^{-\frac{\alpha(I_i(h^f))}{Q_0}h^fx}C_i(0) 
+ \frac{\beta(I_i(h^f))}{\alpha(I_i(h^f))}(1 - e^{-\frac{\alpha(I_i(h^f))}{Q_0}h^f x}).
\end{aligned}
\end{equation}
Since $C_i$ is periodic (i.e.,  $C_i(L)=C_i(0)$),  we get from the previous equation that $C_i(0) = \frac{\beta(I_i(h^f))}{\alpha(I_i(h^f))}$. Inserting this value in~\eqref{eq:Cx}, we find 
\[C_i(x) =C_i^f:=\frac{\beta(I_i(h^f))}{\alpha(I_i(h^f))}, \quad \forall x\in[0, L].\] 
A similar reasoning applied to $p_{C_i}$ gives $p_{C_i}(x) = p_{C_i}^f= \frac{Q_0\gamma(I_i(h^f))}{V N_z\alpha(I_i(h^f))}$, $\forall x\in[0, L]$. It follows that all the terms in the sums of the last equation in~\eqref{eq:optsys} are constant on $[0, L]$. Hence,  there exists a $p_v^f\in\R$ such that the extremal $(C_i^f, v^f, p_{C_i}^f, p_v^f, h^f)$ satisfies the optimality system~\eqref{eq:optsys}. 
\end{proof}

\begin{remark}
The previous theorem shows that the flat topography satisfies the necessary conditions of optimality. One can further explore second-order conditions to check whether the flat topography is a local maximizer. However,  the sign of the eigenvalues of the Hessian operator of the average growth rate Hess$(\bar\mu_{N_z})$ is in general not constant with respect to a flat topography $h^f=V/L$ and is rather difficult to determine (see Appendix~\ref{app:second_order}).

Numerically, we observe that the flat topography is actually optimal in the periodic case for standard values of the parameters (see Subsection~\ref{sec:opt_topo_per}). 
\end{remark}

\begin{remark}\label{rem:non_periodic}
If $C$ is defined by a Cauchy problem and is not assumed to be periodic (i.e., $C(0)$ is not necessarily equal to $C(L)$), then~\eqref{eq:Cx} implies that $C$ may depend on $x$ and the computations in the proof above no longer hold. In other words, the flat topography is not necessarily an optimum in a non-periodic setting,  which is confirmed by our numerical tests (see Subsection~\ref{sec:testC0}).
\end{remark}

\subsection{Non-constant volume problem for maximizing areal productivity}\label{sec:nonconstant}

In the general case, the volume of the system $V$ can also vary, hence can be optimized. We now assume that the water depth is of the form $h + h_0$, where $h\in L^\infty([0, L]; \, \R)$ with $h > -h_0$, $\int_0^L h\, \D x  = 0$, and $h_0 > 0$ so that $V = h_0 L$. Here, $V$ depends only on the parameter $h_0$, as the length $L>0$ is fixed. Moreover, we have $\frac{1}{L}\int_0^L h + h_0\, \D x = \frac{0 + h_0L}{L} = h_0$, meaning that $h_0$ represents the average depth of the system.

On the other hand, when $V$ changes, the biomass concentration $X$ (defined by $\dot X = (\bar \mu - D)X$ with $D$ the dilution rate) also changes. In this case, the light extinction $\varepsilon$ in~\eqref{eq:Beer} can no longer be assumed to be constant. More precisely, we consider here 
\begin{equation}\label{eq:extfun}
\varepsilon(X) := \varepsilon_0 X + \varepsilon_1, 
\end{equation}
where $\varepsilon_0>0$ is the specific light extinction coefficient of the microalgae species and $\varepsilon_1>0$ stands for the background turbidity that summarizes the light absorption and diffusion caused by all non-microalgae components~\cite{Martinez201811}.

To take into account the variation in $X$ with respect to $V$, we also need to adapt our objective function. More precisely, instead of considering the average net growth rate $\bar\mu$, we maximize the areal productivity $\Pi$. Given a biomass concentration $X$, this quantity is defined by 
\begin{equation}\label{eq:Pi}
\Pi := \bar \mu X \frac{V}{S}, 
\end{equation} 
where $\bar \mu$ is the average net growth rate defined in~\eqref{eq:mubar} and $S$ is the ground surface of the raceway system which in our 1D system, actually means $S=L$.

Before stating the associated optimal control problem, we detail the relation between $X$ and $V$. A standard criterion to determine this relation (see~\cite{Masci2010, Grognard2014}) consists in regulating $X$, such that the steady state value of the net growth rate $\mu_s$ at the average depth $h_0$ is 0, i.e.,  
\begin{equation}\label{eq:mus0}
\mu_s\big(I(h_0)\big)=0,  \text{ with }\ \mu_s(I) := -\gamma(I)\frac{\beta(I)}{\alpha(I)} + \zeta(I).
\end{equation}
Using the definitions~\eqref{eq:DefAlphaBeta}, \eqref{eq:DefZetaGamma} for $\alpha$, $\beta$, $\zeta$ and $\gamma$, one can solve~\eqref{eq:mus0} analytically, and find that $I(h_0)$ is one of the two roots, denoted by $I_-$ and $I_+$, of the second-order polynomial equation $k_d\tau R(\sigma I)^2 + (k_r\tau \sigma R-k_r k \sigma)I + k_r R=0$.

In practice, $I_-$, $I_+$ are two real roots with $I_-\leq I_+$, and $\mu_s(I)\geq 0$ on the interval $[I_-, I_+]$. Then, the biomass concentration $X$ in a given volume $V$ is adjusted to get $I(h_0)=I_-$. More precisely, using~\eqref{eq:Beer} with $I(x, z) = I_-$,  we get 
\begin{equation}\label{eq:Xh0}
X(h_0) = \frac{1}{\varepsilon_0}\left(\frac{Y_{\text{opt}}}{h_0} - \varepsilon_1\right), \quad \text{with}\quad Y_{\text{opt}} := \ln\left(\frac{I_s}{I_-}\right).
\end{equation}
Here, $X$ is \Comus{a} function of $h_0$, meaning that we can use the average depth $h_0$ to control both $V$ and $X$ in the non-constant volume case.

\begin{remark}
In bioengineering, the assumption~\eqref{eq:mus0} is usually called the \emph{compensation condition}, which describes the situation where the growth at the bottom compensates exactly for the respiration. We refer to~\cite{Bernard2021optimal} for a detailed analysis.
\end{remark}

We keep using a uniform vertical discretization, as in Section~\ref{sec:objfun}, but now $z_i(0) := \eta(0) - \frac {i-\frac12}{N_z}(h_0+h(0))$, $i=1, \ldots, N_z$. Then the growth rate $\bar\mu_{N_z}$ becomes
\begin{equation}\label{eq:muNzhh0}
\bar \mu_{N_z}(h, h_0) := \sum_{i=1}^{N_z} \int_0^L\frac{\mu\Big(C_i(x), I_i\big(h_0+h(x)\big) \Big)}{h_0 L N_z}(h_0+h(x)) \, \D x.
\end{equation}
Using~\eqref{eq:Xh0} and~\eqref{eq:muNzhh0}, we then derive the semi-discrete areal productivity from~\eqref{eq:Pi}. Note that $V=h_0L$, $X(h_0)$, and $\bar \mu_{N_z}(h, h_0)$ explicitly depend on the average depth $h_0>0$. To treat this parameter, we introduce an additional state variable $y$, such that $y'=0$ and $y=h_0$. This state variable plays the role of $h_0$.

We are now in a position to state the optimal control problem. In the non-constant volume case, we are looking for admissible controls $h\in L^\infty([0, L]; \,\R)$, $h > -y$, and $y>0$ over a fixed length $L > 0$, which maximize the semi-discrete areal productivity. In view of~\eqref{eq:Pi}, the OCP reads as
\begin{equation}\label{eq:ocp2}
\begin{aligned}
\max_{
\begin{aligned}
&h \in L^\infty([0, L]; \mathbb{R})\\
&h > -y, \ y>0
\end{aligned}}  
\Pi_{N_z}(h) := & \sum_{i=1}^{N_z} \int_0^L\frac{\mu\Big(C_i, I_i\big(y+h\big) \Big)}{L N_z}(y+h)X(y) \D x, \\
C'_i = & \frac{\beta(I_i(h+y)) - \alpha(I_i(h+y))C_i}{Q_0}(h+y),\\
C_i(0)= & C_i(L), \qquad \forall i=1, \cdots, N_z,   \\
v'= & h, \\
v(0)= &0,  \ v(L)=0,\\
y' =  &0.
\end{aligned}
\tag{P2}
\end{equation}
Here again, we use formula~\eqref{eq:mu} for $\mu$ and $h$ is the control variable. Moreover, $(C_i, v, y)$ are the state variables, and $X$ is given by~\eqref{eq:Xh0}. The Hamiltonian denoted by $\widetilde{H}$ for the OCP~\eqref{eq:ocp2} is given by
\[\begin{aligned}
\widetilde{H}(C_i, v, y, p_{C_i}, &p_v, p_y, p_0, h) = \sum_{i=1}^{N_z} p_{C_i}\frac{\beta(I_i(h + y)) - \alpha(I_i(h + y))C_i}{Q_0}(h+y)\\
& + p_v h + p_{y}\cdot 0 + p_0 \sum_{i=1}^{N_z} \frac{\mu\Big(C_i, I_i\big(y+h\big) \Big)}{L N_z}(h + y)X(y).
\end{aligned}\]
Here, $(p_{C_i}, p_v, p_y)$ denote the co-states of $(C_i, v, y)$ respectively, and $p_0$ is a real number. Suppose that $h^{\star}\in L^\infty([0, L]; \, R)$, $h^{\star}>-y^{\star}$ is a maximizer, and $(C_i^{\star}, v^{\star}, y^{\star})$ are the corresponding solutions of the problem~\eqref{eq:ocp2}. Using once again the weak maximum principle, there exist absolutely continuous functions $p_{C_i}^{\star}:[0, L]\to \R$, $p_v^{\star}:[0, L]\to \R$, $p_y^{\star}:[0, L]\to \R$ and a real number $p_0^{\star}\leq 0$, such that for almost every $x\in[0, L]$, the extremals $(C_i^{\star}, v^{\star}, y^{\star}, p_{C_i}^{\star}, p_v^{\star}, p_{y}^{\star}, p_0^{\star}, h^{\star})$ satisfy the optimality system
\begin{equation}\label{eq:optsys2}
\begin{aligned}
v' =\ & \frac{\partial \widetilde{H}}{\partial p_v} = h, \quad 
{p_v}' = -\frac{\partial \widetilde{H}}{\partial v} =0, \quad
{y}' = \frac{\partial \widetilde{H}}{\partial p_{y}} =0,\\
{p_{C_i}}' =\ & -\frac{\partial \widetilde{H}}{\partial C_i} 
= p_{C_i} \frac{\alpha\left(I_i(h+y)\right)}{Q_0}(h+y) + p_0\frac{\gamma\left(I_i(h+y)\right)}{LN_z}(h+y)X(y), \\
C_i' =\ & \frac{\partial \widetilde{H}}{\partial p_{C_i}} 
= \frac{\beta\left(I_i(h+y)\right) - \alpha\left(I_i(h+y)\right)C_i}{Q_0}(h+y), \\ 
{p_y}' =\ & -\frac{\partial \widetilde{H}}{\partial y} = 
-\sum_{i=1}^{N_z} p_{C_i}\frac{ \beta'\left(I_i(h+y)\right) - \alpha'\left(I_i(h+y)\right)C_i}{Q_0}(h+y)\partial_{y}I_i(h+y) \\
&- \sum_{i=1}^{N_z} p_{C_i}\frac{\beta\left(I_i(h+y)\right) - \alpha\left(I_i(h+y)\right)C_i}{Q_0}\\
&- p_0\sum_{i=1}^{N_z} \frac{\zeta'\left(I_i(h+y)\right) - \gamma'\left(I_i(h+y)\right)C_i}{LN_z}(h+y)X(y)\partial_{y}I_i(h+y) \\
&- p_0\sum_{i=1}^{N_z} \frac{\zeta\left(I_i(h+y)\right) - \gamma\left(I_i(h+y)\right)C_i}{LN_z}\big(X(y)+(h+y)X'(y)\big) - p_v, \\
0 =\ & \frac{\partial \widetilde{H}}{\partial h} = 
\sum_{i=1}^{N_z} p_{C_i}\frac{ \beta'\left(I_i(h+y)\right) - \alpha'\left(I_i(h+y)\right)C_i}{Q_0}(h+y)\partial_{h}I_i(h+y) \\
&+ \sum_{i=1}^{N_z} p_{C_i}\frac{\beta\left(I_i(h + y)\right) - \alpha\left(I_i(h+y)\right)C_i}{Q_0}\\
&+ p_0\sum_{i=1}^{N_z} \frac{\zeta'\left(I_i(h+y)\right) - \gamma'\left(I_i(h+y)\right)C_i}{LN_z}(h+y)X(y)\partial_{h}I_i(h+y) \\
&+ p_0\sum_{i=1}^{N_z} \frac{\zeta\left(I_i(h+y)\right) - \gamma\left(I_i(h+y)\right)C_i}{LN_z}X(y) + p_v. 
\end{aligned}
\end{equation}

\begin{lemma}\label{lem:normalVnoncst}
The extremals $(C_i^{\star}, v^{\star}, y^{\star}, p_{C_i}^{\star}, p_v^{\star}, p_{y}^{\star}, p_0^{\star}, h^{\star})$ which satisfies~\eqref{eq:optsys2} is normal.
\end{lemma}

\begin{proof}
We follow the same reasoning as in the proof of Lemma~\ref{lem:normalVcst}. Suppose that $p_0^{\star} = 0$, and substitute it into the system~\eqref{eq:optsys2}, the ODE associated with $p_{C_i}^{\star}$ becomes
\[{p_{C_i}}' = p_{C_i} \frac{\alpha\left(I_i(h^{\star} + y^{\star})\right)}{Q_0}(h^{\star} + y^{\star}), 
\quad 
p_{C_i}^{\star}(0)= p_{C_i}^{\star}(L), 
\quad 
\forall i=1, \cdots, N_z. \]
Since $y^{\star} = h_0>0$ and the function $h^{\star} > y^{\star}$, we have $\frac{\alpha\left(I_i(h^{\star} + y^{\star})\right)}{Q_0}(h^{\star} + y^{\star}) > 0$. This implies that $p_{C_i}^{\star} = 0$. Substituting then $p_{C_i}^{\star} = 0$ and $p_0^{\star} = 0$ into the last equation in the system~\eqref{eq:optsys2}, we obtain that $p_{v}^{\star} = 0$, which then implies that ${p_{y}^{\star}}' = 0$. As $p_{y}$ is the co-state associated with the constant $y=h_0$, we have $p_{y}^{\star}(0) = p_{y}^{\star}(L) = 0$, meaning that $p_{y}^{\star}$ also constantly equals 0. Thus, $p_{C_i}^{\star}$, $p_{v}^{\star}$, $p_{y}^{\star}$ are identically 0 on $[0, L]$; which concludes the proof.
\end{proof}

Based on Lemma~\ref{lem:normalVnoncst}, we can normalize the co-states such that $p_0 = -1$. However, unlike Theorem~\ref{thm:flattopo}, the flat topography does not satisfy the optimality system~\eqref{eq:optsys2}.

\begin{theorem}\label{thm:flattopovol}
Given $h_0 > 0$, let $h^f := 0$, $y^f := h_0$, $p_0 = -1$ and assume that $I_s\in (I_-, I_+)$. Then there does not exist a triple $(C_i^{f}, p_{y}^{f}, p_{v}^{f})$ that satisfies the last three equations in the optimality system~\eqref{eq:optsys2}.
\end{theorem}

\begin{proof}
Assuming that there exists such a triple, we start by solving the ODE associated with $C_i^{f}$ in~\eqref{eq:optsys2}. From~\eqref{eq:I_ih}, \eqref{eq:extfun} and~\eqref{eq:Xh0}, we obtain  
\begin{equation}\label{eq:I_ihy}
I_i(h+y)=I_s \exp\Big(- \frac{Y_{\text{opt}}}{y} \frac{i-\frac 12}{N_z} (h+y)\Big),
\end{equation}
where $Y_{\text{opt}}$ is defined in~\eqref{eq:Xh0}. Substituting the values of $h^f$ and $y^f$ into~\eqref{eq:I_ihy}, we find that $I_i(h^f+y^f) = I_i(h_0) = I_s \exp(- Y_{\text{opt}} \frac{i-\frac 12}{N_z})$, which is a constant with respect to $h_0$. A similar analysis to that of the proof of Theorem~\ref{thm:flattopo} shows that $C_i^f=\beta(I_i(h_0))/\alpha(I_i(h_0))$, which is also a constant. Furthermore, differentiating $I_i(h+y)$ with respect to $y$ gives $\partial_{y}I_i(h+y) = I_i(h+y) \cdot \frac{Y_{\text{opt}}}{y^2} \cdot \frac{i-\frac 12}{N_z} h$. Setting $h=h^f$ in this expression, we get $\partial_{y}I_i(h^f+y) = \partial_{y}I_i(0+y) =0$. Substituting all these expressions into the last two equations in~\eqref{eq:optsys2}, we get
\[(p_y^f)' = \frac{X(h_0) + h_0X'(h_0)}{LN_z}\sum_{i=1}^{N_z} \mu_s\big(I_i(h_0)\big) - p_v^f, 
\quad
p_v^f = \frac{X(h_0)}{LN_z}\sum_{i=1}^{N_z}\mu_s\big(I_i(h_0)\big).\]
This implies that $(p_y^f)' = -\frac{Y_{\text{opt}}}{LN_zh_0\varepsilon_0}\sum_{i=1}^{N_z} \mu_s\big(I_i(h_0)\big)$, so that, using~\eqref{eq:Xh0}, we get $X'(h_0) = -\frac{Y_{\text{opt}}}{h_0^2\varepsilon_0}$. Moreover, $I_i(h_0)\in [I_{N_z}(h_0), I_{1}(h_0)] \subset (I_-, I_s) \subset (I_-, I_+)$, hence $\mu_s(I_i(h_0))>0$ for $i\in\{1, \cdots, N_z\}$. We deduce that $(p_y^f)' < 0$. As $p_y^f(0) = p_y^f(L) = 0$, we find a contradiction, which concludes the proof.
\end{proof}

\begin{remark}
Note that the coefficient $h_0$ considered in Theorem~\ref{thm:flattopovol} must satisfy $h_c\leq h_0$ to guarantee that the system remains in a subcritical regime (see Remark~\ref{rem:flu}).
\end{remark}

\section{Numerical Experiments}\label{sec:num}

In this section,  we show some optimal topographies obtained in the various previous frameworks.

\subsection{Numerical Methods}

To solve our optimization problem numerically,  we introduce a supplementary space discretization with respect to $x$. In this way, let us take a space increment $\Delta x$,  set $N_x =[L/\Delta x]$ and $x^{n_x}=n_x\Delta x$ for $n_x=0, \ldots, N_x$. We use Heun's method to compute $(C_i)_{i=1}^{N_z}$ via~\eqref{eq:optsys}. Following a first-discretize-then-optimize strategy,  we get that the co-states $(p_i^C)_{i=1}^{N_z}$ are also computed by a Heun's type scheme. Note that this scheme is still explicit,  since it solves a backward dynamics starting from $p_i(L) = 0$. The optimization is then achieved by a standard gradient method using~\eqref{eq:optsys} and~\eqref{eq:optsys2}, where the stopping criterion involves both the magnitude of the gradient and the constraint $h\geq h_c$,  see Remark~\ref{rem:flu}. The numerical tests are performed by {\sl MATLAB R2020a}~\cite{MATLAB}.

\subsection{Parameter setting}

We now detail the parameters used in our simulations.

\subsubsection{Parameterization}

In our tests, we parameterize $h$ using a truncated Fourier series. More precisely, the water depth reads:
\begin{equation*}
h(x; \boldsymbol{a}) +h_0= h_0 + \sum_{n=1}^N a_n \sin(2n\pi \frac x L),
\end{equation*}
with $\boldsymbol{a} = (a_1, \ldots, a_N)$. This parameterization is motivated by three reasons.
\begin{itemize}
\item  The regularity of the topography is controlled by the order of truncation $N$. As an example, limit situations where $N\rightarrow+\infty$ are not considered in what follows. This framework is consistent with the hydrodynamic regime under consideration, where the solutions of the Saint-Venant equations are smooth.

\item The constraint $h(0; \boldsymbol{a}) = h(L; \boldsymbol{a})$, is preserved, which fits the toric shape of the raceway pond.

\item The water depth has the form $h_0+h$, as assumed in Section~\ref{sec:nonconstant}.
\end{itemize}
From~\eqref{eq:Q0} and~\eqref{eq:zb}, $u$ and $z_b$ also read as functions of $\boldsymbol{a}$. Once the vector $\boldsymbol{a}$ that maximizes $\bar \mu_{N_z}$ is determined, we then find the optimal topography of our system.

\subsubsection{Parameter for the models}

The spatial increment is set to $\Delta x=\SI{0.01}{m}$ so that the convergence of the numerical scheme has been ensured, and we set the raceway length $L=\SI{100}{m}$, the averaged discharge $Q_0 = \SI{0.04}{m^2.s^{-1}}$, the average depth (in the constant volume case) $h_0=h(0;\boldsymbol{a})=\SI{0.4}{m}$ and $z_b(0)=-\SI{0.4}{m}$ to stay in standard ranges for a raceway~\cite{Rayen2019}. The free fall acceleration $g=\SI{9.81}{m.s^{-2}}$. The values of all parameters in Han's model are taken from~\cite{Grenier2020} and given in Table~\ref{tab:parameter}.
\begin{table}
\caption{Parameter values for Han Model }
\label{tab:parameter}
\begin{center}
\begin{tabular}{|c|c|c|}
\hline
$k_r$ & 6.8 $10^{-3}$ & $\si{s^{-1}}$\\
\hline
$k_d$ & 2.99 $10^{-4}$  & -\\
\hline
$\tau$ & 0.25 & $\si{s}$\\
\hline
$\sigma$ & 0.047 & $\si{m^2.\mu mol^{-1}}$\\
\hline
$k$ & 8.7 $10^{-6}$ & -\\
\hline
$R$ & 1.389 $10^{-7}$ & $\si{s^{-1}}$\\
\hline
\end{tabular}
\end{center}
\end{table}

In order to determinate the light extinction $\varepsilon$,  two cases must be considered:
\begin{itemize}
\item constant volume: we assume that only $1\%$ of light can be captured by the cells at the average depth of the raceway,  meaning that $I_-=0.01I_s$, we choose $I_s=\SI{2000}{\mu mol.m^{-2}.s^{-1}}$ which approximates the maximum light intensity, e.g., in summer in the south of France. Then $\varepsilon$ can be computed by $\varepsilon = (1/h_0)\ln(I_s/I_-)$.

\item non-constant volume: in this case, $h_0$ is also a parameter to be optimized. We take from~\cite{Martinez201811} the specific light extinction coefficient of microalgae species $\varepsilon_0=\SI{0.2}{m^2\cdot g}$ and the background turbidity $\varepsilon_1=\SI{10}{m^{-1}}$. 
\end{itemize}

\subsection{Numerical results}

We test the influence of various parameters on optimal topographies. In all of our experiments, we always observe that the obtained topographies satisfy $\min_{x\in[0, L]} h(x; \boldsymbol{a}) > h_c$.

\subsubsection{Influence of vertical discretization}

\begin{figure}
\centering
\includegraphics[scale=0.25]{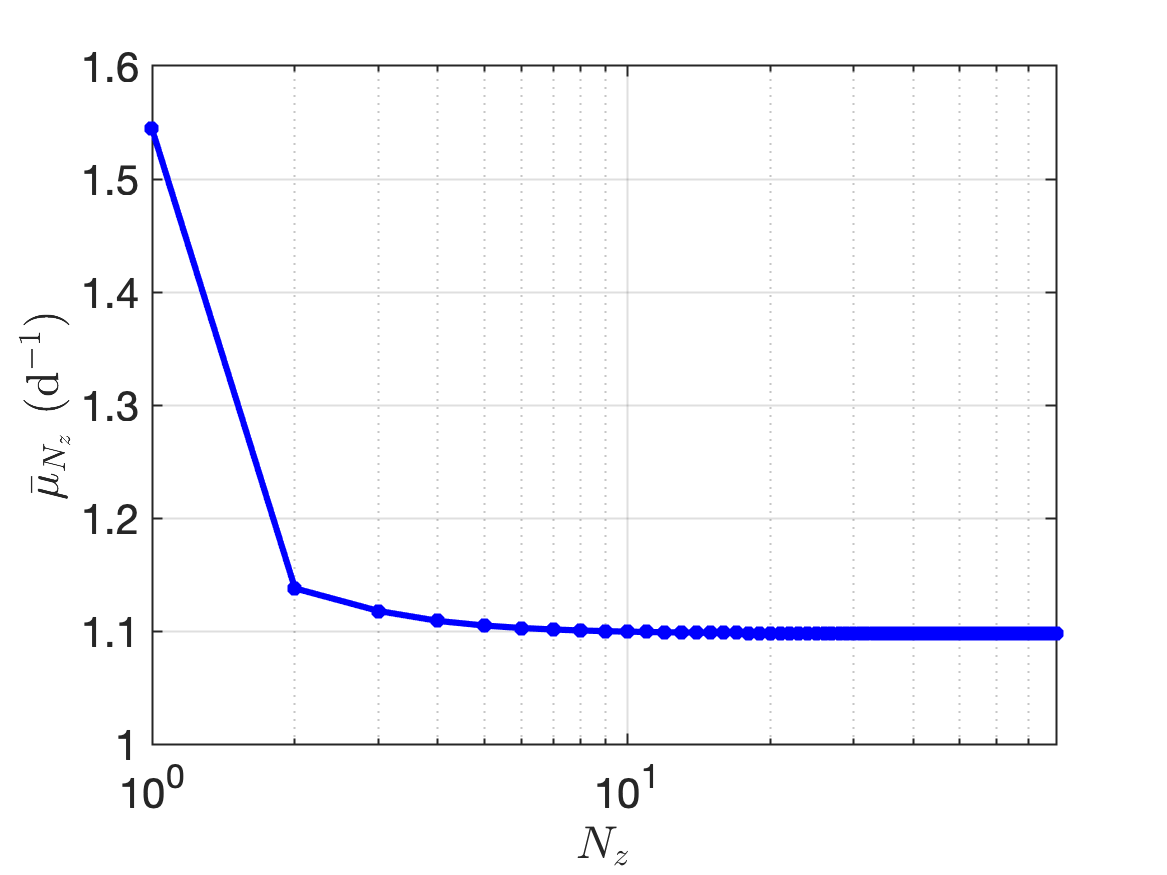}
\caption{Values of the functional $\bar \mu_{N_z}$ for $N_z = [1, 80]$.}
\label{fig:cvNz}
\end{figure}

The first test consists of studying the influence of the vertical discretization parameter $N_z$. We choose $N = 5$, $C_0 = 0.1$ and consider 100 random values $a$. Note that the choice of $a$ should respect the subcritical condition. Let $N_z$ vary from $1$ to $80$, and we compute the average value of $\bar \mu_{N_z}$ for each $N_z$. The results are shown in Fig.~\ref{fig:cvNz}.
We observe numerical convergence when $N_z$ grows, showing the convergence towards the continuous model in space. In view of these results, we take hereafter $N_z = 40$.

\subsubsection{Influence of the initial condition}\label{sec:testC0}

\begin{figure}
\centering
\includegraphics[scale=0.15]{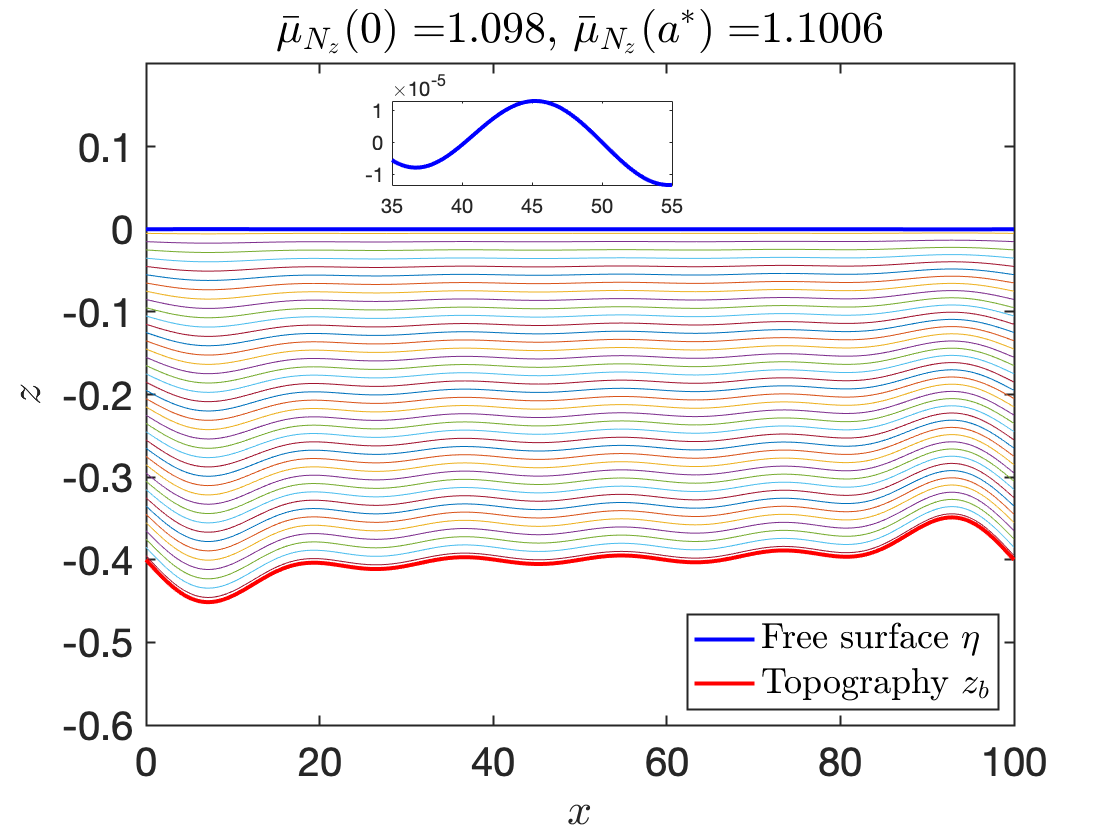}
\includegraphics[scale=0.15]{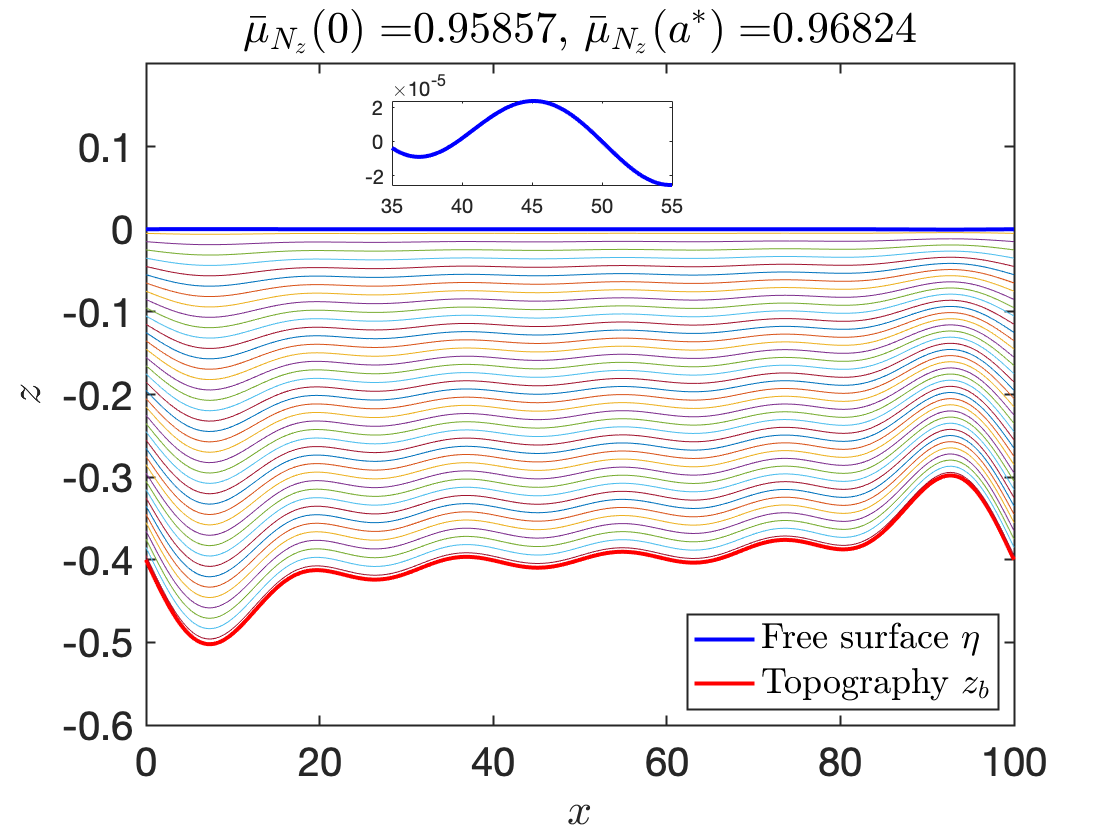}
\caption{Optimal topography for $C_0=0.1$ (left) and $C_0=0.9$ (right). The red thick line represents the topography $z_b$,  the blue thick line represents the free surface $\eta$,  and all the other curves between represent the different trajectories. $\bar \mu_{N_z}(0)$: flat topography,  $\bar \mu_{N_z}(\boldsymbol{a}^*)$: optimal topography.}
\label{fig:C0depend}
\end{figure}

Here, we study the influence of the initial condition $C_0$ on the optimal shape of the raceway pond. We set the numerical tolerance to Tol$ = 10^{-10}$, and consider the order of truncation $N=5$. As for the initial guess, we consider the flat topography, meaning that $\boldsymbol{a}$ is set to 0. We compare the optimal topographies obtained with $C_0=0.1$ and with $C_0=0.9$. The result is shown in Fig.~\ref{fig:C0depend}.
This test confirms Remark~\ref{rem:non_periodic}, since we obtain non-trivial topographies which slightly enhance the algal average growth rate. Moreover, a slight difference between the two optimal topographies is observed. We have observed that this difference remains when the spatial increment $\Delta x$ goes to zero. Although it is difficult to observe in Fig.~\ref{fig:C0depend}, the free surface is not equal to zero, as can be seen for $x\in[35, 55]$.

\subsubsection{Influence of Fourier series truncation}\label{sec:fourier}

The next test is dedicated to the study of the influence of the order of truncation $N$ used to parameterize the water depth $h$. Set $N=[0, 5, 10, 15, 20]$,  $C_0=0.1$ and keep all the other parameters as in the previous section. Table~\ref{tab:N} shows the optimal value of $\bar \mu_{N_z}(\boldsymbol{a}^*)$ and the corresponding maximum eigenvalue of the Hessian $\lambda_{max}(\text{Hess } \bar \mu_{N_z}(\boldsymbol{a}^*))$ for various values of $N$.
\begin{table}
\caption{Behaviour of the objective function for various orders of truncation $N$.}
\begin{center}
\begin{tabular}{|c|c|c|c|c|}
\hline
$N$ & Iter & $\bar \mu_{N_z}(\boldsymbol{a}^*)(\si{d^{-1}})$ & $\log_{10}(\|\nabla \bar \mu_{N_z}(\boldsymbol{a}^*)\|)$ & $\lambda_{max}(\text{Hess } \bar \mu_{N_z}(\boldsymbol{a}^*))$ \\
\hline
0 & 0 & 1.098 & $-$ &  $-$\\
\hline
5 & 16 & 1.1006 & -10.208017 & -6.1400 \\
\hline
10 & 17 & 1.1013 & -10.240885 & -5.9141\\
\hline
15 & 17 & 1.1016 & -10.258798 & -5.9074\\
\hline
20 & 18 & 1.1018 & -10.269413 & -5.9032\\
\hline
\end{tabular}
\end{center}
\label{tab:N}
\end{table}
The result shows a slight increase in the optimal value of $\bar \mu_{N_z}(\boldsymbol{a}^*)$ when $N$ becomes larger. However, the corresponding values of $\bar \mu_{N_z}(\boldsymbol{a}^*)$ remain close to the one associated with a flat topography. Furthermore, the maximum spectrum $\lambda_{max}(\text{Hess } \bar \mu_{N_z}(\boldsymbol{a}^*))$ is always negative, which confirms that local maximizers are obtained.

\subsubsection{Optimal topographies in periodic case}\label{sec:opt_topo_per}

We study the optimal topographies in the constant volume case where the photoinhibition state $C$ is periodic. In our discrete setting,  the Hessian operator is actually of the form $\text{Hess } \bar \mu_{N_z}(h^f)= \lambda Id_{N}$ with $Id_{N}$ the identity matrix of size $N$. We observe that $\lambda <0$,  which confirms that the flat topography is a local maximizer. A precise computation of $\lambda$ together with some remarks about its sign can be found in Appendix~\ref{app:second_order}.

In order to test whether this local maximizer is global, we run the optimization procedure with random admissible topographies. We observe that the procedure always converges to a flat topography (i.e. $\boldsymbol{a}^*=a_f$). This leads us to conjecture that the flat topography corresponds to the global maximum for the average growth rate. For the variable volume case, let us set $N = 5$ (i.e., $\tilde{\boldsymbol{a}}\in\R^6$), and $h_0 = 0.4$ as an initial guess of the average depth. We observe that the optimization stops due to the presence of the physical constraint $h_c$. However, a smaller depth increases the areal productivity, in some cases more than twice the initial areal productivity.

\subsubsection{Simulation with paddle wheel}\label{subsec:twolaps}

In this paragraph, we consider the full raceway pond,  where the mixing induced by the paddle wheel is also considered. More precisely, we simulate several laps with a paddle wheel that mixes up the algae after each lap. The turbulent mixing of the paddle wheel is modeled by a permutation matrix $P$ which rearranges the trajectories in each lap. In our test, $P$ is chosen as an anti-diagonal matrix with entries equal to one. This choice actually corresponds to an optimum and, as shown in~\cite{Bernard02970756}, where other choices are also investigated.

The permutation matrix $P$ corresponds to the permutation $\pi = (1 \ N_z)(2 \ N_z-1)(3 \ N_z-2)\ldots$, where we use the standard notation of cycles in the symmetric group. Note that $\pi$ is of order two. The photoinhibition state $C$ is then set to be $2$-periodic (i.e., $C^1(0) = PC^2(L)$, where $C^1$ and $C^2$ correspond to the photoinhibition state during the first and second lap, respectively). The details of the optimization procedure are given in Appendix~\ref{sec:2laps}.

\begin{figure}
\centering
\includegraphics[scale=0.3]{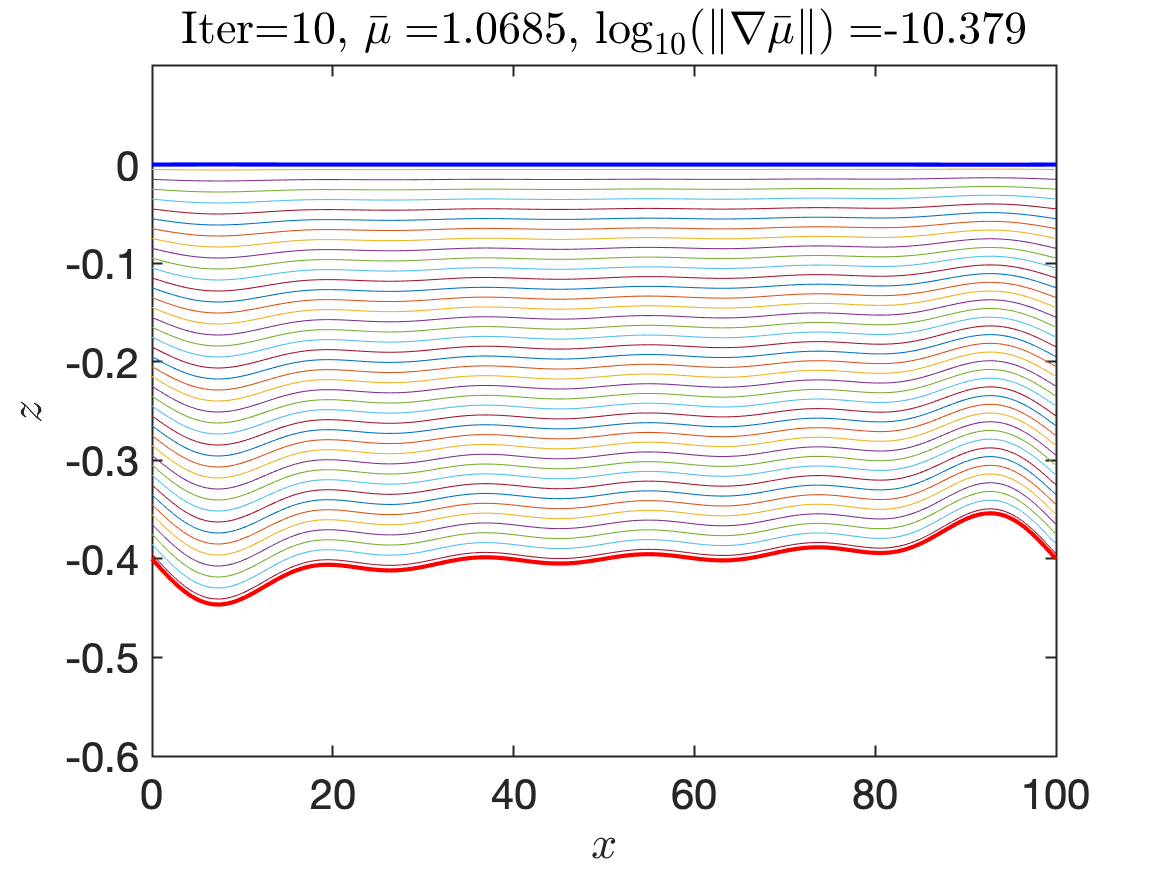}
\includegraphics[scale=0.3]{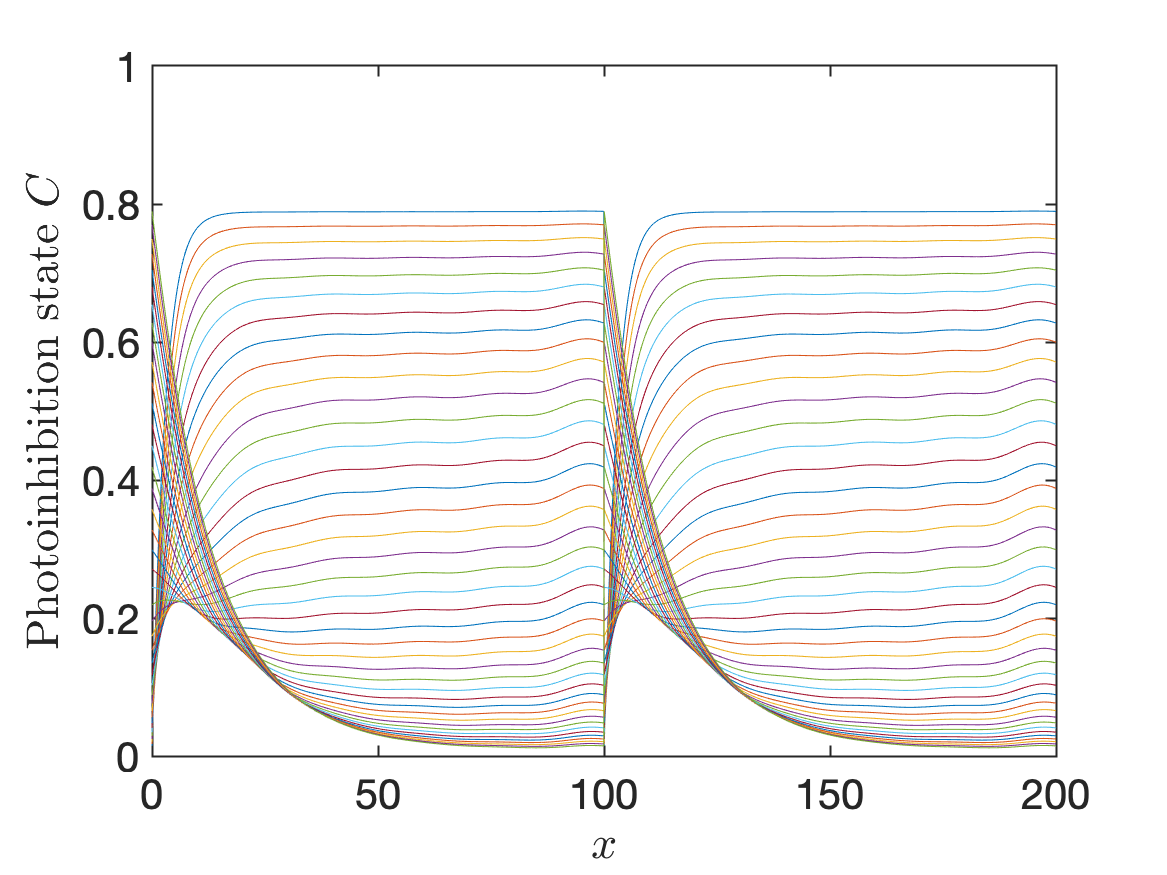}
\caption{Optimal topography (left) and evolution of the photoinhibition state $C$ (right) over two laps.}
\label{fig:2laps}
\end{figure}

We choose a truncation of order $N=5$ in the Fourier series. The initial guess $a$ is set to zero. Fig.~\ref{fig:2laps} presents the shape of the optimal topography and the evolution of the photoinhibition state $C$ over two laps. 
The resulting optimal topography in this case is not flat. However, the increase in the optimal value of the objective function $\bar \mu_{N_z}$ compared to a flat topography with and without permutation are $0.217\%$ and $0.265\%$, respectively, meaning that the increase remains small. On the other hand, we observe that the state $C$ is actually periodic for each lap. This result is actually proved for an arbitrary $P$ in~\cite{Bernard02970756} in the case of a flat topography. This justifies that the optimization strategy only needs to focus on one lap of the raceway (whatever the permutation), and leaves the door open to the optimization of such mixing strategies. We refer to~\cite{Bernard03017414, Bernard03170481} for more details on optimal mixing strategies.

\section{Conclusions and future works}

A flat topography cancels the average algal growth rate gradient when $C$ is assumed to be periodic along the laminar parts of the raceway. This is further confirmed by our numerical tests, in which maximum productivity is obtained for a flat topography. However, considering a more complete framework without periodicity and including a mixing device gives rise to an optimal non-flat topography with a slight gain of the average growth rate. It is not clear whether the difficulty in designing such a pattern could be compensated for by the increase in the process productivity. 

These results may no longer hold if the hydrodynamic regime is turbulent along the entire raceway. In such a case, the increase in the algal productivity may compensate for the higher energetic cost of mixing. However, without the laminar assumption, the problem becomes challenging,  and much work remains to be done in this direction.

{\bf Acknowledgments.} We are very grateful to Emmanuel Tr\'elat (Laboratoire Jacques-Louis Lions, Sorbonne Universit\'e, Paris) for the helpful discussions we have had on this topic. We thank the editors and anonymous referees for their valuable and constructive comments, which greatly improved the quality of this paper.

\bibliographystyle{siamplain}        
\bibliography{auto}

\appendix
\section{Two-lap system with a paddle-wheel}\label{sec:2laps}

Denote by $P$ the permutation matrix associated with $\pi=(1 \ N_z)(2 \ N_z-1)(3 \ N_z-2)\ldots$ (see Section~\ref{subsec:twolaps}), i.e., 1 as entries on the anti-diagonal and by $C^1$ (resp. $C^2$) the photoinhibition state for the first (resp. second) lap of the raceway. We then assume that the state $C$ is 2-periodic, meaning that $C^1(0)=PC^2(L)$. From~\eqref{eq:muNz},  we define the objective function by 
\[\frac12\sum_{j=1}^2\bar \mu^j_{N_z}(h) = \frac12\sum_{j=1}^2\sum_{i=1}^{N_z}\int_0^L\frac{\mu\Big(C_i^j(x), I_i\big(h(x)\big) \Big)}{V N_z}h \, \D x.\]
For a fixed volume $V>0$ and a discharge $Q_0>0$, the associated OCP reads:
\begin{equation}\label{eq:cond_2laps}
\begin{aligned}
\max_{h\in L^\infty(0, L; \; \R), \ h>0} \frac12\sum_{j=1}^2 \bar \mu_{N_z}(h)=\ & \frac12\sum_{j=1}^2\sum_{i=1}^{N_z}\int_0^L\frac{\mu\Big(C_i^j(x), I_i\big(h(x)\big) \Big)}{V N_z}h \, \D x,\\
{C_i^j}' =\ & \frac{\beta\left(I_i(h)\right) - \alpha\left(I_i(h)\right)C_i^j}{Q_0}h, \\
C^1(L)=\ & PC^2(0),  \quad C^1(0) = PC^2(L),  \\
v'=\ & h, \\
v(0)=\ &0,  \ v(L)=V.
\end{aligned}
\end{equation}
Denote by $H$ the Hamiltonian associated with this problem, which reads
\[\begin{aligned}
H(C_i^j, v, p_{C_i}^j, p_v, p_0, h) = &\sum_{j=1}^2\sum_{i=1}^{N_z} p_{C_i}^{j}\Big( \frac{\beta\left(I_i(h)\right) -\alpha(I_i(h))C_i^j}{Q_0}h\Big) + p_v h\\
& + p_0 \frac12\sum_{j=1}^2\sum_{i=1}^{N_z}\frac{\zeta\left(I_i(h)\right) - \gamma\left(I_i(h)\right)C_i^j}{V N_z}h, 
\end{aligned}\]
where $p_{C_i}^{j}$, $p_v$ are the co-states of $C_i$, $v$, and $p_0$ is a real number. A similar analysis to that of Section~\ref{sec:cst_vol} gives a similar optimality system as~\eqref{eq:optsys}, in which $p_{C_i}^{j}$ satisfies the conditions $p_{C}^1(L)= Pp_{C}^2(0)$ and $p_{C}^2(L)= Pp_{C}^1(0)$.

\section{Second order conditions}\label{app:second_order}
Consider the second-order condition under the truncated Fourier parameterization. Since Fourier modes $(\sin(2n\pi\frac xL))_{n\in\mathbb{N}}$ are orthogonal, a direct computation gives Hess $\bar \mu_{N_z}(h^f)= \lambda Id_{N}$ with
\[\begin{aligned}
\lambda = \frac1{Q_0}\sum_{i=1}^{N_z} &2 p_{C_i} \big(\beta'(I_i(h))-\alpha'(I_i(h)) C_i\big) I'_i(h) 
+ p_{C_i}\big(\beta'(I_i(h))-\alpha'(I_i(h)) C_i\big)I''_i(h) h\\
& + p_{C_i}\big(\beta''(I_i(h))-\alpha''(I_i(h)) C_i\big){I'_i(h)}^2h\\
+ \frac{p_0}{V N_z}&\sum_{i=1}^{N_z} 2\big(\zeta'(I_i(h)) - \gamma'(I_i(h))C_i\big)I'_i(h) 
+ \big(\zeta'(I_i(h)) - \gamma'(I_i(h))C_i\big)I''_i(h)h \\
& + \big(\zeta''(I_i(h)) - \gamma''(I_i(h))C_i\big){I'_i(h)}^2h.
\end{aligned}\]
Using the definitions~\eqref{eq:DefAlphaBeta} and~\eqref{eq:DefZetaGamma}, we get $\alpha(I)=\beta(I)+k_r$ and $\zeta(I)=\gamma(I)-R$. As $\alpha'(I)=\beta'(I)$ and $\zeta'(I)=\gamma'(I)$, one gets
\begin{equation}\label{eq:lambda}
\begin{aligned}
\lambda = \sum_{i=1}^{N_z}&(1-C_i)\Big[\frac{p_{C_i}}{Q_0}\Big(2 \beta'(I_i(h)) I'_i(h) 
+ \beta'(I_i(h))I''_i(h) h + \beta''(I_i(h)){I'_i(h)}^2h\Big) \\
&+ \frac{p_0}{V N_z}\Big(2\gamma'(I_i(h))I'_i(h) + \gamma'(I_i(h))I''_i(h)h 
+ \gamma''(I_i(h)){I'_i(h)}^2h\Big)\Big].
\end{aligned}
\end{equation}
Furthermore, one can differentiate the closed forms of $I(h)$,  $\beta(I)$ and $\gamma(I)$ to have 
\[\begin{aligned}
&I'_i(h) = -\varepsilon \frac{i-\frac 12}{N_z} I_i(h),  \quad I''_i(h)=(\varepsilon\frac{i-\frac 12}{N_z})^2 I_i(h),  \\ 
&\beta''(I)=\frac2{(\tau\sigma I+1)(\tau\sigma I+2)I}\beta'(I), \quad \gamma''(I)=-\frac{2\sigma\tau}{\tau\sigma I+1}\gamma'(I).
\end{aligned}\]
Inserting these analytical forms into~\eqref{eq:lambda} gives
\[\begin{aligned}
\lambda = \sum_{i=1}^{N_z}&(1-C_i)\varepsilon \frac{i-\frac 12}{N_z} I_i(h)\Big[\frac{p_{C_i}\beta'(I_i(h))}{Q_0}(h\varepsilon \frac{i-\frac 12}{N_z} 
+ \frac{2h\varepsilon \frac{i-\frac 12}{N_z}}{(\tau\sigma I_i(h)+1)(\tau\sigma I_i(h)+2)} - 2) \\
&+\frac{p_0\gamma'(I_i(h))}{V N_z}\big(h\varepsilon \frac{i-\frac 12}{N_z} - \frac{2\sigma \tau h\varepsilon \frac{i-\frac 12}{N_z} I_i(h)}{\tau\sigma I_i(h)+1} - 2\big)\Big].
\end{aligned}\]
Considering now the case $h=h^f=V/L$, one gets 
\[\begin{aligned}
&1-C_i^f = \frac{k_r}{\alpha(I_i(h^f))}>0,  
\quad p_{C_i}^f = p_0\frac{Q_0\gamma(I_i(h^f))}{V N_z\alpha(I_i(h^f))} <0,  \\
&\beta'(I)=\frac{k_d\tau\sigma^2 I(I\sigma\tau + 2)}{(I\sigma\tau + 1)^2}>0, 
\quad \gamma'(I)=\frac{k\sigma}{(I\sigma\tau + 1)^2}>0.
\end{aligned}\]
Hence, in the limit case, the sign in the big bracket becomes positive when $h$ goes to 0 and the flat topography is no longer a local maximizer for small values of $h$ in this case. Under the assumption that the hydrodynamics is subcritical, then $\lambda<0$ in practice as shown in Section~\ref{sec:fourier} and in Section~\ref{sec:opt_topo_per}.
\end{document}